\documentclass[preprint,letterpaper,11pt,oneside]{elsarticle}
\usepackage{amssymb,amsmath,amsthm}
\usepackage{latexsym,epsfig,multicol,anysize,graphicx}
\usepackage{latexsym,amsfonts,amssymb,amscd}
\usepackage{fancybox,shadow,color,graphicx,psfrag,float}
\usepackage{graphicx,epsfig,tikz}
\usepackage{geometry}
\usepackage{verbatim}
\usepackage{xcolor}

\usepackage[colorlinks=true,linkcolor=red,
citecolor=red, urlcolor=green]{hyperref}

\newcommand{\floor}[1]{{\left\lfloor #1 \right\rfloor}}
\newcommand{\ceil}[1]{{\left\lceil #1 \right\rceil}}

\def\f{\mathbb F}

\def\a{\alpha}
\def\b{\beta}

\def\ii{\infty}

\def\N0{\mathbb{N}_0}

\newcommand{\cha}{\operatorname{char}}

\newcommand{\Sup}{\operatorname{Supp}}
\newcommand{\res}{\operatorname{res}}

\newtheorem{teo}{Theorem}[section]

\newtheorem{pro}[teo]{Proposition}
\newtheorem{lem}[teo]{Lemma}
\newtheorem{cor}[teo]{Corollary}

\newtheorem{example}{Example}

\usepackage{etoolbox}
\makeatletter
\newcommand{\changeoperator}[1]{%
  \csletcs{#1@saved}{#1@}%
  \csdef{#1@}{\changed@operator{#1}}%
}
\newcommand{\changed@operator}[1]{%
  \mathop{%
    \mathchoice{\textstyle\csuse{#1@saved}}
               {\csuse{#1@saved}}
               {\csuse{#1@saved}}
               {\csuse{#1@saved}}%
  }%
}
\makeatother

\changeoperator{sum}
\changeoperator{prod}

\usepackage{scalerel}

\makeatletter
\def\ps@pprintTitle{%
  \let\@oddhead\@empty
  \let\@evenhead\@empty
  \let\@oddfoot\@empty
  \let\@evenfoot\@oddfoot
}
\makeatother

\begin{document}
\begin{frontmatter}
\title{Bases for Riemann--Roch spaces of linearized function fields with applications to generalized algebraic geometry codes}
\author{Horacio Navarro}
\ead{horacio.navarro@correounivalle.edu.co}
\date{}
\address{Departamento de Matem\'aticas, Universidad del Valle, Cali-Colombia.}

\begin{abstract}

In this paper, we determine explicit bases for Riemann--Roch spaces of linearized function fields, and we give a lower bound for the minimum distance of generalized algebraic geometry codes. As a consequence, we construct generalized algebraic geometry codes with good parameters.
\end{abstract}
\end{frontmatter}





\section{Introduction}

In 1999, Xing, Niederreiter, and Lam \cite{xing99} presented a construction of linear codes based on algebraic curves over finite fields that generalizes the algebraic geometry codes (AG codes) introduced two decades before by Goppa. These codes are called generalized algebraic geometry codes (GAG codes), and the main feature of this construction is that it is possible to use places of arbitrary degree. In 2002, Heydtmman \cite{heydtmann2002generalized} studied how to define and represent the dual code of a GAG code, and she gave a decoding algorithm for this kind of codes. Later, Dorfer and Maharaj gave a short proof of the duality result presented by Heydtmman, and they showed generator and parity check matrices for GAG codes. To calculate those matrices as well as to determine the parameters of GAG codes, Riemann--Roch spaces play a fundamental role, as in Goppa's construction. The paper is organized as follows. 

In Section 2, we recall the definition of GAG codes, their properties, and present a generalization of the lower bound for the minimum distance given by Picone in \cite{picone2013}. In Section 3, we show the main properties of linearized function fields and calculate explicit bases of the Riemann--Roch spaces. Analogous results for function fields of Kummer type were proved by Garz\'on and Navarro in \cite{navarro2022}. In Section 4, we show examples of GAG codes with good parameters, using the results obtained in the previous sections.

\section{Generalized algebraic geometry codes }

The following notation will be used throughout this correspondence. Let $F/\mathbb{F}_q$ be an algebraic function field with genus $g$.  The Riemann--Roch space corresponding to the divisor $G$ is  defined as \[\mathcal{L}(G)=\{x\in F: (x)\geq -G\}\cup \{0\}, \]
where $(x)$ is the principal divisor of $x$. $\mathcal{L}(G)$ is a finite-dimensional vector space over $\mathbb{F}_q$, and its dimension is denoted by $\ell(G)$. 

The differentials space of the divisor $G$ is defined as
\[\Omega(G)=\{\eta \in \Omega: (\eta)\geq G\}\cup \{0\}, \]
where $\Omega$ denotes the differentials space of $F$. $\Omega(G)$ is also a finite-dimensional vector space over $\mathbb{F}_q$, and its dimension, called the index of specialty of $G$,  is  \[i(G)=\ell(G)-\deg G+g-1.\]
For a canonical divisor $W$ of $F$, $i(G)$ is also equal to $\ell(W-G)$.  Let $P_1,\ldots,P_s$ be distinct places of $F$ of degrees $k_1,\ldots,k_s$, respectively, and $G$ be a divisor with support disjoint from the support of the divisor $D:=\sum_{i=1}^sP_i$. Let $C_i$ be an $[n_i,k_i, d_i]$ linear code  over $\mathbb{F}_q$, and  let $\pi_i:\mathbb{F}_{q^{k_i}}\rightarrow C_i \subseteq \mathbb{F}_q^{n_i}$ be a fixed  $\mathbb{F}_q$-linear  isomorphim mapping $\mathbb{F}_{q^{k_i}}$ onto  $C_i$ for any $i=1,2,\ldots,s$.  Let $n=\sum_{i=1}^sn_i$. The generalized algebraic geometry codes (GAG codes)  $C_\mathcal{L}(P_1,\ldots,P_s;G;C_1,\ldots,C_s)$ and $C_\Omega(P_1,\ldots,P_s;G;C_1,\ldots,C_s)$ are defined as the images of the $\mathbb{F}_q$-linear maps

\[
\alpha_{\mathcal{L}}:  \mathcal{L}(G) \longrightarrow  \mathbb{F}_q^n, \quad 
x  \mapsto  (\pi_1(x(P_1)),\pi_2(x(P_2)), \ldots, \pi_s(x(P_s)))
\]

  \[\alpha_{\Omega}:  \Omega(G-D) \longrightarrow  \mathbb{F}_q^n, \quad 
\omega  \mapsto  (\pi_1(\res_{P_1}(\omega)), \pi_2(\res_{P_2}(\omega)),\ldots, \pi_s(\res_{P_s}(\omega)))
.\]

If  we take rational places  and for each $1\leq i\leq s$ we choose $C_i=[1,1,1]$  the trivial lineal code over $\mathbb{F}_q$ and $\pi_i$ the identity map on $\mathbb{F}_q$, then the GAG codes reduce to AG codes.

Now, we recall some results concerning to the parameters of GAG codes.

\begin{teo} $C_\mathcal{L}(P_1,\ldots,P_n;G;C_1,\ldots,C_s)$ is an $[n,k,d]$ code with parameters 
\[
k=\ell(G)-\ell(G-D)\quad \text{ and } \quad
d\geq \sum_{i=1}^s d_i-\max \left\{\sum_{i \in S}d_i:S \in X \right\},
\]
where \[X=\left\{S\subseteq \{1,\ldots, s \} :\sum_{i\in S}k_i\leq \deg G \right\}. \]
Furthermore, if  $\deg G < \deg D$, then $k=\ell(G)\geq \deg G +1-g$, and the equality holds if $2g-2<\deg G < \deg D$.
\end{teo}

\begin{teo}\label{teo1GAG}$C_{\Omega}(P_1,\ldots,P_s;G;C_1,\ldots,C_s)$ is an $[n,k,d]$ code with parameters
\[k=i(G-D)-i(G)=\ell(G-D)-\ell(G)+\deg D \quad \text{ and } \]
\[d \geq \sum_{i=1}^s d_i-\max \left\{\sum_{i \in S}d_i:S \in X \right\},\]
where
\[X=\left\{S\subseteq \{1,\ldots, s \} :\sum_{i \in S} k_i\leq \deg D -\deg G +2g-2 \right\}. \]

Furthermore, if  $\deg G > 2g-2$ then \[k=\ell(G-D)-\deg G-1+g+\deg D\geq \deg D+g -\deg G-1\] and the equality holds if $2g-2<\deg G < \deg D$.
\end{teo}

We state the lower bound for the minimum distance of a GAG code given by Picone.

\begin{teo}(\cite[Proposition 2.9]{picone2013})\label{distGAGpicone} Let $F/\mathbb{F}_q$ be an algebraic function field of genus $g$ and $P_1,\ldots,P_s$ distinct places of $F$ of degrees $k_1,\ldots,k_s$ respectively, and $C_i$ be an $[n_i,k_i, d_i]$ linear code  over $\mathbb{F}_q$ for any $i=1,2,\ldots,s$.
Let  $A$, $B$, $C$, and $Z$ be divisors of $F$ with support disjoint from the support of $D=\sum_{i=1}^sP_i$. Suppose   $\mathcal{L}(A) = \mathcal{L}(A-Z)$ and $\mathcal{L}(B) = \mathcal{L}(B + Z)= \mathcal{L}(C)$. If $G = A + B$ and $d_i\geq k_i$ for $1\leq i \leq s$, then the minimum distance $d$ of the code $C_{\Omega}(P_1,\ldots,P_s;G;C_1,\ldots,C_s)$  satisfies
\[d \geq \deg G-(2g-2)+\deg Z -i(A)+i(G-C). \]

\end{teo}

We remove some hypotheses from the result above to obtain a more general lower bound.

Here, it is necessary to introduce the notion of the greatest common divisor of two divisors. Given $D_1$ and $D_2$, two divisors of a function field $F/\f_q$, the \textit{greatest common divisor} of $D_1$  and  $D_2$ is defined as \[\gcd(D_1,D_2):=\sum_{Q\in \mathbb{P}_F} \min\{\nu_Q(D_1), \nu_Q(D_2)\}Q. \]

\begin{teo}\label{distAGgen1}  Let $F/\mathbb{F}_q$ be an algebraic function field of genus $g$ and $P_1,\ldots,P_s$ distinct places of $F$ of degrees $k_1,\ldots,k_s$ respectively, and $C_i$ be an $[n_i,k_i, d_i]$ linear code  over $\mathbb{F}_q$ for any $i=1,2,\ldots,s$.
Let  $A$, $B$, $C$, and $Z$ be divisors of $F$ with support disjoint from the support of $D=\sum_{i=1}^sP_i$. Suppose   $\mathcal{L}(A) = \mathcal{L}(A-Z)$, $\mathcal{L}(B) = \mathcal{L}(B + Z)\subseteq \mathcal{L}(C)$. If $G = A + B$, then the minimum distance $d$ of the code $C_{\Omega}(P_1,\ldots,P_s;G;C_1,\ldots,C_s)$ satisfy
\[d \geq \sum_{i=1}^s d_i-\max \left\{\sum_{i \in S}d_i:S \in X\right\}, \]
where \[X=\left\{S\subseteq \{1,\ldots, s \} :\sum_{i \in S} k_i\leq \deg D -\deg G -\deg Z +2g-2-i(A)+i(G-C) \right\}. \]

\end{teo}
\begin{proof}
Let $\omega \in \Omega(G-D)$  such that $\alpha_{\Omega}(\omega) $ is a non-zero codeword of minimum weight.  Let $S$ be the set of all $i$ such that $\res_{P_i}(\omega)=0$ and let $T$ be the complement of $S$ in $\{1,2,\ldots,s\}$.
Then 
\begin{equation}\label{equcota1}
d\geq \sum_{i \in T} d_i=\sum_{i=1}^s d_i-\sum_{i \in S} d_i.
\end{equation}
Consider the divisor  $D':=\sum_{i \in T} P_i$. Since  $\alpha_{\Omega}(\omega) $ is zero for $i\in S$, then $\omega\in \Omega(G-D')$, that is, $W:=(\omega)\geq G-D'$. Thus, $E:=W-G+D'$ is an effective divisor such that $\Sup(E)\cap\Sup(D')=\emptyset $ and we have $W=G-D'+E$. Taking degrees of both sides, we get that 
\[2g-2=\deg G-\deg D' +\deg E =\deg G -\sum_{i \in T} \deg P_i +\deg E,\]
then   
\[\sum_{i \in T} \deg  P_i -\deg G +2g-2-\deg E\geq 0,\]
Summing $\sum_{i \in S} \deg P_i =\sum_{i \in S} k_i$ in both sides  of the inequality above, we have
\begin{equation*}
\deg D-\deg G + 2g-2 -\deg E  \geq  \sum_{i \in S} k_i,
\end{equation*}
hence
\begin{equation}\label{equcota2}
\deg D-\deg G + 2g-2 -\deg Z -i(A)+i(G-C)  \geq  \sum_{i \in S} k_i,
\end{equation}
if we assume $\deg E \geq \deg Z +i(A)-i(G-C) $. Finally, from \eqref{equcota1} and  \eqref{equcota2} the desired result follows. It remains to prove that $\deg(E)\geq \deg(Z)+i(A)-i(G-C)$. Notice that 
\begin{align}\label{equcota6}
\deg E &=\ell(A+E)-\ell(A)+i(A)-i(A+E) \nonumber \\ 
&=\ell(A+E)-\ell(A-Z) +i(A)-i(A+E) \nonumber  \\
&\geq \ell(A+E)-\ell(A+E-Z)+i(A)-i(A+E).
\end{align} 
By the Riemann--Roch Theorem, we get \[\ell(A+E)=\deg(A+E)+1-g+\ell(W-(A+E))\]
and 
\[\ell(A+E-Z)=\deg(A+E-Z)+1-g+\ell(W-(A+E-Z)),\]
hence,
\begin{align}\label{equcota7}
\ell(A+E)-\ell(A+E-Z)&=\deg Z+\ell(W-(A+E))-\ell(W-(A+E-Z)) \nonumber \\
&=\deg Z+\ell(B-D')-\ell(B+Z-D') 
\end{align}
Since $D'$ is effective and $B,C$ and  $Z$ have support disjoint from the support of $D'$, we get  
\[\gcd(B+Z-D',B)\leq B-D' \text{ and } \gcd(B-D',C)\leq C-D'.\]
As \[\mathcal{L}(B+Z-D')\subseteq \mathcal{L}(B+Z)=\mathcal{L}(B),\]
then
\[\mathcal{L}(B+Z-D')= \mathcal{L}(B+Z-D')\cap\mathcal{L}(B)=\mathcal{L}(\gcd(B+Z-D',B))\subseteq\mathcal{L}(B-D')\]
and thus  
\begin{equation}\label{equcota8}
\ell(B+Z-D')\leq \ell(B-D').
\end{equation}

In a similar way,  note that 
\[\mathcal{L}(B-D')\subseteq \mathcal{L}(B)\subseteq \mathcal{L}(C),\]
implies 
\[\mathcal{L}(B-D')= \mathcal{L}(B-D')\cap\mathcal{L}(C)=\mathcal{L}(\gcd(B-D',C))\subseteq \mathcal{L}(C-D')\]
and  hence
\begin{equation}\label{equcota9}
\ell(B-D')\leq \ell(C-D').
\end{equation}

Finally, we have 
\begin{align}\label{equcota10}
i(A+E)&=\ell(W-(A+E))=\ell(B-D')\leq  \ell(C-D')\,\, (\text{by } \eqref{equcota9} ) \nonumber  \\
&\leq \ell(C-D'+E)= i(W-C+D'-E)=i(G-C). 
\end{align}

Combining this with \eqref{equcota6}--\eqref{equcota8} yields the desired result.

\end{proof}

\begin{cor}\label{distAGgen}Let $F/\mathbb{F}_q$ be an algebraic function field of genus $g$ and $P_1,\ldots,P_s$ distinct places of $F$ of degrees $k_1,\ldots,k_s$ respectively, and $C_i$ be a $[n_i,k_i, d_i]$ linear code  over $\mathbb{F}_q$ for any $i=1,2,\ldots,s$. Let  $A$, $B$, and $Z$ be divisors of $F$ with support disjoint from the support of $D=\sum_{i=1}^sP_i$. Suppose   $\mathcal{L}(A) = \mathcal{L}(A-Z)$ and $\mathcal{L}(B) = \mathcal{L}(B + Z)$. If $G = A + B$, then the minimum distance $d$ of the code $C_{\Omega}(P_1,\ldots,P_s;G;C_1,\ldots,C_s)$ satisfies
\[d \geq \sum_{i=1}^s d_i-\max \left\{\sum_{i \in S}d_i:S \in X\right\}, \]
where \[X=\left\{S\subseteq \{1,\ldots, s \} :\sum_{i \in S} k_i\leq \deg D -\deg G -\deg Z +2g-2 \right\}. \]
\end{cor}

\begin{proof}
This follows immediately from Proposition  \ref{distAGgen1},  taking $C=B$.
\end{proof}

We use this result in Section \ref{examples} to construct examples of GAG codes with good parameters.

\section{Bases for Riemann--Roch spaces of linearized function fields}

We start this section by defining linearized function fields and showing some important properties that will be used frequently.

A \textit{linearized function field} is a function field $F=\mathbb{F}_{q^n}(x,y)$ defined by the equation 
\begin{equation}\label{linearized}
L(y)=h(x):=\frac{f(x)}{g(x)},
\end{equation}
with $L(y)=\sum_{i=0}^{r}a_iy^{q^i}\in \mathbb{F}_{q^n}[y]$ a  
linearized  polynomial,  $a_r, a_0\neq 0$, and  $q^r$  roots in  $\mathbb{F}_{q^n}$.  We assume  throughout that the factorization into irreducibles of $f$ and $g$ is $f=\prod_{j=1}^mq_j^{m_j},$ and $ g=\prod_{i=1}^sp_i^{n_i}$ respectively,   and the positive integers  $n_i$, $d:=\deg f -\deg g>0$  are coprime to $\cha \mathbb{F}_q$.

\begin{pro}\label{propo1} With the notations above, we have the following:

\begin{enumerate}[i.]
\item For $1\leq i\leq s$, let  $P_i$ (resp. $P_\infty$) be the zero of $p_i$ (resp. the pole of $x$) in $\mathbb{F}_{q^n}(x)$. The places $P_1,\ldots, P_s,$ and $ P_\infty$ are  totally ramified places in $F/\mathbb{F}_{q^n}(x)$.

\item\label{int. basis1} The functions $1,y,\ldots, y^{q^r-1}$ form an integral basis  at any place $P$ of $\mathbb{F}_{q^n}(x)$ different from $P_1,\ldots, P_s,$ and $P_\infty$. Moreover, $P_1,\ldots, P_s,$ and $P_\infty$ are the only ramified places in $F/\mathbb{F}_{q^n}(x)$. 

\item\label{propo1.3} For $1\leq i\leq s$, let  $Q_i$ (resp. $Q_\infty$) be the only zero of $p_i$ (resp. the only pole of $x$) in $F$. Then the principal divisor of $p_i$  in $F$  is
\[(p_i)=q^rQ_i-q^r\deg p_i Q_{\infty},\]

\item \label{propo1.4} For $1\leq j\leq m$, let $S_j$ be the only zero of $q_j$ in $\mathbb{F}_{q^n}(x)$. Then $S_j$ splits completely in $F$. Moreover, the principal divisor of $y$ in $F$  is
\[(y)=\sum_{j=1}^m m_jR_j-\sum_{i=1}^s n_iQ_i-(\deg f-\deg g)Q_{\infty},\]
where  $R_j$ is the only zero of $y$ in $F$ lying over $S_j$, for each $1\leq j \leq m$. 
\end{enumerate}

\end{pro}

\begin{proof}\begin{enumerate}[i.]
\item Let $Q$ be a place lying over $P_\infty$ and let $e$ be the index of ramification of $Q|P_\infty$. Since $\nu_{P_\infty}(h(x))=-(\deg f-\deg g)=-d<0$ we have $\nu_{Q}(L(y))=q^r\nu_{Q}(y)$ and 
\begin{equation}\label{ecuuuu9}
q^r\nu_{Q}(y)=\nu_{Q}(L(y))=e\nu_{P_\infty}(h(x))=-e\cdot d.
\end{equation}
By the equation \eqref{ecuuuu9} and as $\gcd(d,q)=1$, then $P_\infty$ is totally ramified and $\nu_{Q}(y)=-d.$ In the same way as before, we have that $P_i$ is totally ramified and  $\nu_{Q}(y)=-n_i,$ with  $1\leq i \leq s.$ 
\item Let $P$ be any place of $\mathbb{F}_{q^n}(x)$ different from $P_1,\ldots, P_s$, and $P_\infty$, and let $Q$ be a place of $F$ lying over $P$. It is clear that the minimum polynomial of $y$ over $\mathbb{F}_{q^n}(x)$ satisfies
  \[\varphi(t)=L(t)-h(x)=t^{q^r}+a_{r-1}t^{q^{r-1}} +\cdots +a_1t^{q}+a_0t-h(x)\in \mathcal{O}_P[T],\]
and $\nu_Q(\varphi'(y))=\nu_Q(a_0)=0$, then by [\cite{stichbook}, Corollary 3.5.11] $\{1,y,\ldots, y^{q^r-1}\}$ is an integral basis of $F$ at $P$ and $P$ is not ramified in $F/\mathbb{F}_{q^n}(x)$.

\item Since the principal divisor of $p_i$ in $\mathbb{F}_{q^n}(x) $ is $(p_i)=P_i-\deg p_i Q_\infty$ and the places $P_i$ and $P_{\infty}$ are totally ramified, then \[(p_i)=\sum_{Q\in\mathbb{P}_F}\nu_Q(p_i)Q=\nu_{Q_{i}}(p_i)Q_{i}+\nu_{Q_{\infty}}(p_i)Q_{\infty}=q^rQ_i-q^r\deg p_i Q_{\infty}.\]

\item The place $S_j$ is a zero of $h(x)$ because
 it is the zero of $q_j$ in $\mathbb{F}_q(x)$, then $\varphi(t) \in \mathcal{O}_{S_j}[t]$. The reduction modulus $S_j$ of $\varphi(t)$ is  the polynomial
\[\overline{\varphi}(t)=t^{q^r}+a_{r-1}t^{q^{r-1}} +\cdots +a_1t^q+a_0t,\]
which has all its roots in $\mathbb{F}_{q^n}$, hence  there exist $q^r$ places lying over $S_j$ by [\cite{stichbook}, Theorem 3.3.7]. One of them is a zero of $y$ lying over $S_j$, and we will call it $R_j$.  It is clear that $\nu_{R_j}(y)=m_j$. All of the above allows us  to conclude that \[(y)=\sum_{Q\in\mathbb{P}_F}\nu_Q(p_i)Q=\sum_{j=1}^m m_jR_j-\sum_{i=1}^s n_iQ_i-(\deg f-\deg g)Q_{\infty}.\]
\end{enumerate}
   
\end{proof}

We can now formulate our main results.

\begin{teo}\label{generatorset}Consider the divisor $G:=\sum_{i=1}^s a_iQ_i + b_{\infty}Q_{\infty}$ on the function field defined by the equation \eqref{linearized} where $a_i, b_{\infty} \in \mathbb{Z}$ for $1\leq i \leq s$. Then 
\[S:=\left\{y^kx^{e_{\infty,k}}\prod_{i=1}^sp_i^{e_{i,k}}:
\begin{array}{c}
0\leq e_{\infty,k}, \,  e_{i,k}\in \mathbb{Z},  -a_i\leq -kn_i+q^re_{i,k}, \\
\\
  k(\deg f-\deg g)+q^r( e_{\infty,k}+\sum_{i=1}^s e_{i,k}\deg p_i)\leq  b_\infty, \\
  \\
   \text{ for all } k, \,  0\leq k \leq q^r-1.
\end{array} 
\right\}\]

is a generator set of $\mathcal{L}(G).$

\end{teo}

\begin{proof} 
We first show that $S\subseteq \mathcal{L}(G)$. By Proposition \ref{propo1}.\ref{propo1.3},  \ref{propo1.4} and since the divisor of $x$ is $(x)=(x)_0-q^rQ_{\infty}$, we have that the divisor of the function $y^kx^{e_{\infty,k}}\prod_{i=1}^s p_i^{e_{i,k}} \in S$ satisfies the inequality \[\left(y^kx^{e_{\infty,k}}\prod_{i=1}^s p_i^{e_{i,k}}\right)\geq \sum_{i=1}^s(q^re_{i,k}-kn_i)Q_i-\left(kd+q^r( e_{\infty,k}+\sum_{i=1}^s e_{i,k}\deg p_i)\right)Q_\infty,\]
with $d=\deg f-\deg g$ and for any $0\leq k\leq q^r-1$. Since $q^re_{i,k}-kn_i\geq -a_i$ and $kd+q^r( e_{\infty,k}+\sum_{i=1}^s e_{i,k}\deg p_i)\leq b_\infty$, then we conclude that $S\subseteq \mathcal{L}(G)$ and  the $\mathbb{F}_{q^n}$-linear span of $S$ is a subset of $\mathcal{L}(G)$.

Next, we prove that every element of $\mathcal{L}(G)$ can be expressed as a linear combination of elements from $S$ with coefficients from $\mathbb{F}_{q^n}$. Let $z\in \mathcal{L}(G)$ be  a non-zero element. Then the places $Q_i$ and $Q_{\infty}$ are the only possible poles of $z$. Since $\{1,y,\ldots,y^{q^r-1}\}$ is a basis of $F/\mathbb{F}_{q^n}(x)$,  there exist $z_i\in \mathbb{F}_{q^n}(x) $ such that \[z=z_0+z_1y+\cdots+z_{q^r-1}y^{q^r-1},\]
and  the only possible poles in $\mathbb{F}_{q^n}(x)$ of the $z_k$ are the places $P_1,\ldots,P_s$ and $P_{\infty }$. Indeed, suppose that 
$P\not \in \{P_1,\ldots,P_s,P_{\infty}\}$  is a pole of $z_k$ for some $0\leq k\leq q^r-1$. Note that for  any place $Q$ lying over $P$, we have  $Q \not\in \Sup(G)$, which gives $\nu_Q(z)\geq 0$  and   $z\in \bigcap_{Q|P}\mathcal{O}_Q$. 
By [\cite{stichbook}, Corollary 3.3.5] and  Proposition \ref{propo1}.\ref{int. basis1}  the integral closure $\mathcal{O}'_P$ of $\mathcal{O}_P$ in $F$ is 
\[\mathcal{O}'_P=\bigcap_{Q|P}\mathcal{O}_Q=\sum_{k=0}^{q^r-1}\mathcal{O}_Py^k.\]
Hence, $z\in \mathcal{O}'_P$, so 
each  $z_k$ must belong to $\mathcal{O}_P$. This  contradiction implies that the $z_k$ are of the form 
\[z_k=f_k\prod_{i=1}^{s}p_i^{e_{i,k}}\]
where $e_{i,k}$ are integers, $f_k$ is  a polynomial in $\mathbb{F}_{q^n}[x]$, and $p_i$ does not divide $f_k$ for any $1\leq i \leq s$. Thus, $z_k$  is a $\mathbb{F}_{q^n}$-linear combination of functions 
\[A_{k,m}=x^m\prod_{i=1}^{s}p_i^{e_{i,k}}\]
for $m=0,1,\ldots, \deg f_k$. In order to prove the theorem, it is sufficient to show that for $0\leq k\leq q^r-1$ and $m=0,1,\ldots, \deg f_k$, the functions $y^kA_{k,m}$ belong to $S$. Let  $0\leq j \leq s$, $Q:=Q_j$, and $P:=P_j$.  Since $P$ is the place defined by the irreducible polynomial $p_j\in \mathbb{F}_{q^n}[x]$, we have  
\begin{equation*}
\nu_P(z_k)=\nu_P(f_k)+\nu_P\left(\prod_{i=1}^s p_i^{e_{i,k}}\right)=e_{j,k}.
\end{equation*}
On the other hand, by Proposition \ref{propo1}.\ref{propo1.4} $\nu_Q(y)=-n_j$, then 
\[\nu_Q(z_ky^k)=q^r\nu_P(z_k)+k\nu_Q(y)=q^r\nu_P(z_k)-kn_j=q^re_{j,k}-kn_j.\]

Now, for $0\leq k\leq q^r-1$, we have that  $\nu_Q(z_ky^k)$ are distinct because they are distinct module $q^r$. Thus,  
\begin{equation*}
q^re_{j,k}-kn_j\geq \min_{0\leq k\leq q^r-1}\{\nu_Q(z_ky^k)\}=\nu_Q(z)\geq -a_j
\end{equation*}
since $z\in \mathcal{L}(G)$.
Finally,  for $0\leq k\leq q^r-1$ and $m=0,1,\ldots, \deg f_k$, we obtain 
\begin{eqnarray*}
\nu_Q(y^kA_{k,m})&=&\nu_Q(y^kx^m\prod_{i=1}^{s}p_i^{e_{i,k}})
=   k\nu_{Q}(y)+m\nu_{Q}(x)+\sum_{i=1}^s{e_{i,k}}\nu_{Q}(p_i) \\
&\geq & k\nu_{Q}(y)+e_{j,k}\nu_{Q}(p_j)
 =    q^re_{j,k} -kn_j \geq -a_j.
\end{eqnarray*}
It remains to establish 
$\nu_{Q_\infty}(y^kA_{k,m})\geq -b_\infty.$  We proceed in a similar way as before. Recalling $z_k=f_k\prod_{i=1}^sp_i^{e_{i,k}}$, we get 
\begin{equation}\label{ecuacion 2.6}\nu_{P_\infty}(z_k)=\nu_{P_\infty}(f_k)+\sum_{i=1}^s{e_{i,k}}\nu_{P_\infty}(p_i) =-\deg f_k-\sum_{i=1}^s e_{i,k}\deg p_i.
\end{equation}

Now, for $0\leq k\leq q^r-1$, the integers $q^r\nu_{P_\infty}(z_k)-kd$ are distinct because they  are distinct module $q^r$. Then, 
\begin{equation}\label{ecuacion 2.7}
q^r\nu_{P_\infty}(z_k)-kd\geq \min_{0\leq k \leq q^r-1}\{\nu_{Q_\infty}(z_ky^k)\}=\nu_{Q_\infty}(z)\geq -b_\infty,
\end{equation}
because $z\in \mathcal{L}(G)$.
On the other hand,  for $m=0,1,\ldots, \deg f_k$ we have
\begin{eqnarray}\label{ecuacion 2.8}
\nu_{Q_\infty}(y^kA_{k,m})&=& k\nu_{Q_\infty}(y)+m\nu_{Q_\infty}(x)+\sum_{i=1}^s{e_{i,k}}\nu_{Q_\infty}(p_i)\nonumber  \\
&=&-kd-q^r\left(m+\sum_{i=1}^s e_{i,k}\deg p_i\right).
\end{eqnarray}

Then, combining \eqref{ecuacion 2.6}--\eqref{ecuacion 2.8}, we obtain
\begin{eqnarray*}
-\nu_{Q_\infty}(y^kA_{k,m})&=&kd+q^r\left(m+\sum_{i=1}^s e_{i,k}\deg p_i\right) \\
& \leq & kd+q^r\left(\deg f_k+\sum_{i=1}^s e_{i,k}\deg p_i\right) \\
&\leq & b_\infty.
\end{eqnarray*}
\end{proof}

\begin{lem}[\cite{maharaj2005riemann}, Lemma 3.5]\label{lema clave}

Let $F/K$ be a function field. Let $G$ be a divisor of $F$, and let $P$ be a rational place of $F$. Let $V=\{\nu_P(z): z\in \mathcal{L}(G)\setminus\{0\}\}$. For each $i \in V$, choose $u_i\in \mathcal{L}(G)$  such that $\nu_P(u_i)=i.$ Then the set $B=\{u_i:i \in V\}$ is a basis for $\mathcal{L}(G)$.
\end{lem}

\begin{teo}\label{maintheo}Consider the divisor $G:=\sum_{i=1}^s a_iQ_i + b_\infty Q_\infty$ on the function field defined by the equation \eqref{linearized} where $a_i, b_\infty \in \mathbb{Z}$ for $1\leq i \leq s$. Then the dimension of $\mathcal{L}(G)$ over $\mathbb{F}_{q^n}$ is

\[\ell(G)=\sum_{k=0}^{q^r-1}\max\left\{\sum_{i=1}^{s}\floor{\frac{a_i-kn_i}{q^r}}\deg p_i+\floor{\frac{b_\infty-kd}{q^r}}+1,0\right\},\]
with $d=\deg f-\deg g$.

\end{teo}

\begin{proof}

For $0\leq k\leq q^r-1$ set 
\[V_k:=\left\{-kd-q^r\left(e+\sum_{i=1}^s e_i\deg p_i\right):
\begin{array}{c}  
 e_i\in \mathbb{Z},\, -a_i\leq -kn_i+q^re_i \text{ for all }  1 \leq i \leq s,\\ 
 0 \leq e,\, kd+q^r(e+\sum_{i=1}^s e_i\deg p_i)\leq  b_\infty.
\end{array}
\right\} \]
and let $V := \bigcup_{k=0}^{q^r-1}V_k$. The proof will be  divided into several steps.
\vspace{0.5cm}

\textit{Step 1.} $V=\{\nu_{Q_\infty}(z):z\in \mathcal{L}(G)\setminus \{0\}\}$.
\vspace{0.5cm}

\textit{Proof of step 1.}
Let a non-zero element $z\in \mathcal{L}(G)$.  By Theorem \ref{generatorset}, we have  
 \[z=\sum_{k=0}^{q^r-1}\alpha_k\left( y^k x^{e_{\infty, k}}\prod_{i=1}^sp_i^{e_{i,k}} \right)=\sum_{k=0}^{q^r-1}\left( \alpha_k x^{e_{\infty, k}}\prod_{i=1}^sp_i^{e_{i,k}} \right)y^k:=\sum_{k=0}^{q^r-1}z_k y^k \] where  $\alpha_k \in \mathbb{F}_{q^r}$, 
 $e_{\infty,k}$ is a nonnegative integer, each $ e_{i,k}$ is an integer, for all $1 \leq i \leq s$,   $-a_i\leq -kn_i+q^re_{i,k}$ and $kd+q^r(e_{\infty,k}+\sum_{i=1}^s e_{i,k}\deg p_i)\leq  b_\infty $.  Thus,  for any $0 \leq k \leq q^r-1$, we obtain    
\[\nu_{Q_\infty}(z)=\nu_{Q_\infty}(z_ky^k)=q^r \nu_{P_\infty}(z_k)+  k\nu_{Q_\infty}(y)= -q^r\left(e_{\infty,k}+\sum_{i=1}^s e_{i,k}\deg p_i\right)-kd \in V_k.  \]
Conversely, let $m\in V$. Then $m\in V_j$ for some $0\leq j \leq q^r-1$. Hence, 
  \[m=-kd-q^r\left(e+\sum_{i=1}^s e_i\deg p_i\right)\] 
where $0\leq e$, $e_i\in \mathbb{Z}$, $-a_i\leq -kn_i+q^re_i$, and $kd+q^r\left(e+\sum_{i=1}^s e_i\deg p_i\right)\leq  b_\infty$. Finally, if we choose \[z=y^kx^e\prod_{i=1}^sp_i^{e_i}\] we have $ \nu_{Q_\ii}(z)=m$.  
\vspace{0.5cm}

According to Lemma \ref{lema clave}, it follows that $\ell(G)=|V|$. Therefore, we proceed to count the number of elements in $V$. This will be the second step of the proof.

\textit{Step 2.} Fix $k$, $0\leq k \leq q^r-1.$ Then $-kd -cq^r\in V$ if and only if 

\[-\sum_{i=1}^{s}\floor{\frac{a_i-kn_i}{q^r}}\deg p_i\leq c \leq \frac{b_\infty -kd}{q^r}\cdot\]

\textit{Proof of step 2.} Firstly, we will prove that the sets $V_k$ are mutually disjoint. Let $N\in V_k\cap V_j$, then $-kd-q^ra=N=-jd-q^rb$, where $a,b \in \mathbb{Z}$,  so 
\[kd\equiv jd \mod q^r.\]
Since $(d, q)=1$, then $k \equiv j \mod q^r $, thus $k=j$. This implies that \[\ell(G)=|V|=\sum_{k=0}^{q^r-1}|V_k|.\]
It remains to determine $|V_k|$.
Now fix $k$,  $0\leq k \leq q^r-1$, and set $m:=-kd-q^rc$ where $c\in \mathbb{Z}.$  Then $m\in V$ if and only if $m=-kd -q^rc\in V_k$. This holds if and only if 

\[-kd -q^rc=-kd-q^r\left(e+\sum_{i=1}^s e_i\deg p_i\right), \]
for some integers $e_i$ and $e\geq 0$ with
\[kd+q^r(e+\sum_{i=1}^s e_i\deg p_i)\leq  b_\infty\] and  \[ -a_i\leq -kn_i+q^re_i\quad \text{  for any }\quad   1\leq i \leq s.\] These inequalities are equivalent to 
\[c=e+\sum_{i=1}^s e_i\deg p_i\]
and 
\[ e+\sum_{i=1}^s e_i\deg p_i\leq \frac{b_\infty-kd}{q^r}\]

and \[e_i\geq  \ceil{-\frac{a_i-kn_i}{q^r}}=-\floor{\frac{a_i-kn_i}{q^r}}\]  for any  $1\leq i \leq s$. Therefore, the integers $e_i$ and $e\geq0$  exist if and only if 

\[ -\sum_{i=1}^{s}\floor{\frac{a_i-kn_i}{q^r}}\deg p_i \leq c\leq \frac{b_\infty-kd}{q^r}.\]
Now it follows that $|V_k|$ is the number of integers in the interval \[ \left[-\sum_{i=1}^{s}\floor{\frac{a_i-kn_i}{p^r}}\deg p_i, \floor{\frac{b_\infty-kd}{p^r}}\right],\] 
i.e.,
\[|V_k|:=\max\left\{\sum_{i=1}^{s}\floor{\frac{a_i-kn_i}{q^r}}\deg p_i+\floor{\frac{b_\infty-kd}{q^r}}+1,0\right\}.\]  
This implies that \[\ell(G)=\sum_{k=0}^{q^r-1}\max\left\{\sum_{i=1}^{s}\floor{\frac{a_i-kn_i}{q^r}}\deg p_i+\floor{\frac{b_\infty-kd}{q^r}}+1,0\right\}\] 
\end{proof}

\begin{cor}\label{basesteo} Consider the divisor $G:=\sum_{i=1}^s a_iQ_i + b_\infty Q_\infty$ on the function field defined by the equation \eqref{linearized} where $a_i, b_\infty \in \mathbb{Z}$ for $1\leq i \leq s$. Then 

\[\bigcup_{0\leq k \leq q^r-1}\mathcal{B}_k\]

where each $\mathcal{B}_k$ is defined as

\[\left\{y^kx^{e_{\infty, k}}\prod_{i=1}^sp_i^{e_{i,k}}: 
 0\leq e_{_{\infty,k}}, 
-\sum_{i=1}^s \floor{\frac{a_i-kn_i}{q^r}}\deg p_i \leq e_{_{\infty, k}}+ \sum_{i=1}^s e_{_{i,k}}\deg p_i\leq  \floor{ \frac{b_\infty-kd}{q^r}} \right\}\] is a basis for $\mathcal{L}(G)$ as a vector space over $\mathbb{F}_{q^n}$. 
\end{cor}

\begin{proof}
This follows immediately from Theorem \ref{generatorset} and the proof of Theorem \ref{maintheo}.
\end{proof}

\begin{cor}[\cite{stichbook}, Proposition 4.1 (h)] Consider a function field $F=\mathbb{F}_{q^n}(x,y)$ with \[y^q+\mu y=g(x)\in \mathbb{F}_{q^n}[x], \] 
$0\neq \mu \in \mathbb{F}_{q^n}$. Assume that  $\deg g:=m$ is relatively prime to $q$ and   $\{\beta: \beta^q+\mu \beta=0\}\subseteq \mathbb{F}_{q^n}$. Then a basis for the Riemann-Roch space  $\mathcal{L}(rQ_{\infty})$ is 
\[ \bigcup_{\scaleto{k=0}{5pt}}^{\scaleto{q-1}{5pt}}\mathcal{B}_k=\{y^k x^j:0\leq j , \quad 0\leq k \leq q-1,  \quad jq+km\leq r \}. \]
\end{cor}

Now, we recall the notion of the floor of a divisor. The \textit{floor} of a divisor  $G$ with $\ell(G)>0$, denoted $\floor{G}$,   is the minimum degree divisor such that $\mathcal{L}(G) = \mathcal{L}(\floor{G})$. The next result states that from a generator set of a Riemann--Roch space $\mathcal{L}(G)$, we can calculate $\floor{G}$. 

\begin{lem}\label{floordiv}[\cite{maharaj2005riemann}, Theorem 2.6] Let $G$ be a divisor of $F/ \mathbb{F}_q$, and let $b_1, \ldots , b_t \in \mathcal{L}(G)$ be a spanning set
for $\mathcal{L}(G)$. Then,
$\floor{G}=-\gcd((b_i ):i = 1, \ldots , t).$
\end{lem}

In order to conclude this section, we compute the floor of certain divisors of linearized function fields.

\begin{teo}\label{floorteo} Let $G:=\sum_{i=1}^s a_iQ_i + b_\infty Q_\infty$ be a divisor of the function field defined by the equation \eqref{linearized} where $a_i, b_\infty \in \mathbb{Z}$ for $1\leq i \leq s$. The floor of $G$ is given by \[\floor{G}:=\sum_{i=1}^s \a_iQ_i + \b_\infty Q_\infty,\]
with

\[\a_i:=\max_{\scaleto{0\leq k \leq q^r-1}{5pt}}\left\{kn_i+q^r\floor{\frac{a_i-kn_i}{q^r}}: \sum_{i=1}^{s}\floor{\frac{a_i-kn_i}{q^r}}\deg p_i+\floor{\frac{b_\infty-kd}{q^r}}\geq 0  \right\}\]
and 
\[\b_{\infty}:=\max_{\scaleto{0\leq k \leq q^r-1}{5pt}}\left\{kd+ q^r\floor{\frac{b_\infty-kd}{q^r}}:\sum_{i=1}^{s}\floor{\frac{a_i-kn_i}{q^r}}\deg p_i+\floor{\frac{b_\infty-kd}{q^r}}\geq 0  \right\}\]
and  $d=\deg f-\deg g$.
\end{teo}

\begin{proof} By  Corollary \ref{basesteo}
, $\mathcal{B}=\bigcup_{k=0}^{q^r-1}\mathcal{B}_k$ is a basis for $\mathcal{L}(G)$, where each $\mathcal{B}_k$ is the set
\[
\left\{y^kx^{e_{\scaleto{\infty, k}{4pt}}}\prod_{i=1}^sp_i^{e_{i,k}}: 
0\leq e_{\scaleto{\infty, k}{4.5pt}}, 
-\sum_{i=1}^s \floor{\frac{a_i-kn_i}{q^r}}\deg p_i \leq e_{\scaleto{\infty, k}{4.5pt}}+ \sum_{i=1}^s e_{_{i,k}}\deg p_i\leq \floor{\frac{b_{_\infty}-kd}{q^r}} 
\right\}
.\]
Suppose $\mathcal{B}_k$, is non-empty and consider $z=y^kx^e\prod_{i=1}^sp_i^{e_i} \in \mathcal{B}_k$. By Proposition \ref{propo1}.\ref{propo1.3},  \ref{propo1.4}, and since the divisor of $x$ is $(x)=(x)_0-q^rQ_{\infty}$ and $k\sum_{j=1}^m m_jR_j +(x)_0$ is an effective divisor, we have
\begin{eqnarray*}
(z)&=&k\sum_{j=1}^m m_jR_j+e(x)_0+ 
\sum_{i=1}^s(-kn_i+e_iq^r)Q_i-(kd+ q^r(e +\sum_{i=1}^se_i \deg p_i)) Q_{\infty}\\
&\geq & \sum_{i=1}^s\left(-kn_i -q^r\floor{\frac{a_i-kn_i}{q^r}}\right)Q_i-\left(kd+ q^r\floor{\frac{b_\infty-kd}{q^r}}\right)Q_{\infty}.\\
\end{eqnarray*}
Now, take $z=y^k\prod_{i=1}^sp_i^{e_i}$ with $e_i=-\floor{\frac{a_i-kn_i}{q^r}}$. Then 
\[\nu_{Q_i}(z)=-kn_i+e_iq^r=-kn_i-q^r\floor{\frac{a_i-kn_i}{q^r}}.\] 

Observe that $\mathcal{B}_k$ is non-empty if and only if $\sum_{i=1}^{s}\floor{\frac{a_i-kn_i}{q^r}}\deg p_i+\floor{\frac{b_\infty-kd}{q^r}}\geq 0$. Therefore,

\begin{eqnarray*}
\alpha_i&=&-\min_{0\leq k \leq q^r-1}\left\{\nu_{Q_i}(z):\sum_{i=1}^{s}\floor{\frac{a_i-kn_i}{q^r}}\deg p_i+\floor{\frac{b_\infty-kd}{q^r}}\geq 0\right\}\\
&&\\
&=&\max_{0\leq k \leq q^r-1}\left\{kn_i+q^r\floor{\frac{a_i-kn_i}{q^r}}: \sum_{i=1}^{s}\floor{\frac{a_i-kn_i}{q^r}}\deg p_i+\floor{\frac{b_\infty-kd}{q^r}}\geq 0 \right\} .\\
\end{eqnarray*}
It remains to calculate $\beta_{\infty}$.  Similarly, choosing $z=y^kx^{e}\prod_{i=1}^sp_i^{e_i}$ with $e_i=-\floor{\frac{a_i-kn_i}{q^r}}$ and \[e=\floor{\frac{b_\infty-kd}{q^r}}+\sum_{i=1}^{s}\floor{\frac{a_i-kn_i}{q^r}}\deg p_i=\floor{\frac{b_\infty-kd}{q^r}}-\sum_{i=1}^{s}e_i\deg p_i,\] we can see that
\[\nu_{Q_\ii}(z)=-kd- q^r(e +\sum_{i=1}^se_i \deg p_i)=-kd-q^r\floor{\frac{b_\infty-kd}{q^r}}.\]
Thus,
\begin{eqnarray*}
\beta_\infty&=&-\min_{0\leq k \leq q^r-1}\left\{\nu_{Q_\ii}(z):\sum_{i=1}^{s}\floor{\frac{a_i-kn_i}{q^r}}\deg p_i+\floor{\frac{b_\infty-kd}{q^r}}\geq 0 \right\}\\
&=&\\
&=&\max_{0\leq k \leq q^r-1}\left\{kd+q^r\floor{\frac{b_\infty-kd}{q^r}}:\sum_{i=1}^{s}\floor{\frac{a_i-kn_i}{q^r}}\deg p_i+\floor{\frac{b_\infty-kd}{q^r}}\geq 0 \right\},\\
\end{eqnarray*}

and by Lemma \ref{floordiv}, we obtain
 $\floor{G}=\sum_{i=1}^s \a_iQ_i + \b_\infty Q_\infty.$
\end{proof}

\section{Examples}\label{examples}
In this section, we show examples of GAG codes with good parameters. For this purpose, we choose divisors $A$, $B$, and $R$ such that $Z=A-\floor{A}$, $Z=R-\floor{R}$, and $B=\floor{R}$. Then \[\ell(A-Z)=\ell(A)\quad \text{ and } \quad\ell(B+Z)=\ell(B) \] In this way, $A,B,Z$ and $G=A+B$ satisfy the hypotheses of the Corollary \ref{distAGgen}, so we obtain lower bounds for minimum distance. We use Theorem \ref{teo1GAG} to calculate the dimension of the codes, Theorem \ref{maintheo} for the dimension of the Riemann--Roch spaces involved, and Theorem \ref{floorteo} for the floor of the divisors considered. All examples presented are improvements on MinT's tables \cite{MinT}. In examples 1–3, we consider the function field $F=\mathbb{F}_{49}(x,y)$ defined by the equation
\[y^7+y=\frac{(x^2+1)^2}{x^2}.\]
In \cite{gupta2023reciprocal}, it was proven that $F/\mathbb{F}_{49}$ has genus $12$ and $170$ rational places. Also, it can be shown that it has many places of degree two and three. The only pole of $x$ and $y$ and the only pole of $y$ over the zero of $x$ in $\mathbb{F}_{49}(x)$ will be denoted by $P_{\infty}$ and $P_0$, respectively.
\begin{example} Let $A=5P_\infty+18P_0$, $B=4P_\infty+18P_0$, $G=A+B=9P_\infty+36P_0$, and $Z=P_{\infty}$, then $\ell(G)=34$.
\begin{enumerate} \item Taking $109\leq s \leq 168$ and $[n_i,k_i,d_i]=[1,1,1]$ for $1\leq i\leq s$, we obtain linear codes with parameters $[s, s-34, 24]$.
\item Taking $116\leq s \leq 168$ and $[n_i,k_i,d_i]=[1,1,1]$ for $1\leq i\leq s$ and $[n_i,k_i,d_i]=[2,2,1]$ for $i=s+1$, we obtain linear codes with parameters $[s+2, s-32, 23]$.
\item Taking $139\leq s \leq 168$ and $[n_i,k_i,d_i]=[1,1,1]$ for $1\leq i\leq s$ and $[n_i,k_i,d_i]=[2,2,1]$ for $i=s+1, s+2$, we obtain linear codes with parameters $[s+4, s-30, 22]$.

\item Taking $112\leq s \leq 168$ and $[n_i,k_i,d_i]=[1,1,1]$ for $1\leq i\leq s$ and $[n_i,k_i,d_i]=[3,2,2]$ for $i=s+1$, we obtain linear codes with parameters $[s+3, s-32, 24]$.
\item Taking $123\leq s \leq 168$ and $[n_i,k_i,d_i]=[1,1,1]$ for $1\leq i\leq s$ and $[n_i,k_i,d_i]=[3,2,2]$ for $i=s+1,s+2$, we obtain linear codes with parameters $[s+6, s-30, 24]$.

\item Taking $123\leq s \leq 168$ and $[n_i,k_i,d_i]=[1,1,1]$ for $1\leq i\leq s$ and $[n_i,k_i,d_i]=[3,2,2]$ for $i=s+1,s+2$, we obtain linear codes with parameters $[s+6, s-30, 24]$.
\item Taking $163\leq s \leq 168$ and $[n_i,k_i,d_i]=[1,1,1]$ for $1\leq i\leq s$ and $[n_i,k_i,d_i]=[5,4,2]$ for $i=s+1$, we obtain linear codes with parameters $[s+5, s-30, 22]$.
\end{enumerate} \end{example}
\begin{example} Let $A=19P_\infty+4P_0$, $B=18P_\infty+3P_0$, $G=A+B=37P_\infty+7P_0$, and $Z=P_\infty$. Then $\ell(G)=33$. \begin{enumerate} \item Taking $107\leq s \leq 168$ and $[n_i,k_i,d_i]=[1,1,1]$ for $1\leq i\leq s$, we obtain linear codes with parameters $[s,s-33,23]$. 
\item Taking $127\leq s \leq 168$ and $[n_i,k_i,d_i]=[1,1,1]$ for $1\leq i\leq s$ and $[n_i,k_i,d_i]=[2,2,1]$ for $i=s+1$, we obtain linear codes with parameters $[s+2, s-31, 22]$.
\item Taking $150\leq s \leq 168$ and $[n_i,k_i,d_i]=[1,1,1]$ for $1\leq i\leq s$ and $[n_i,k_i,d_i]=[2,2,1]$ for $i=s+1, s+2$, we obtain linear codes with parameters $[s+4, s-29, 21]$.

\item Taking $127\leq s \leq 168$ and $[n_i,k_i,d_i]=[1,1,1]$ for $1\leq i\leq s$ and $[n_i,k_i,d_i]=[3,2,2]$ for $i=s+1,s+2$, we obtain linear codes with parameters $[s+6, s-29, 23]$.

\end{enumerate}

\end{example}

\begin{example}
Let $A=19P_\infty$, $B=18P_\infty+4P_0$, $G=A+B=37P_\infty+4P_0$ and $Z=P_\infty$. Then $\ell(G)=30$. 

\begin{enumerate}
\item Taking $129\leq s \leq 168$ and $[n_i,k_i,d_i]=[1,1,1]$ for $1\leq i\leq s$, we obtain linear codes with parameters $[s,s-30,20]$.

\item Taking $152\leq s \leq 168$ and $[n_i,k_i,d_i]=[1,1,1]$ for $1\leq i\leq s$  and  $[n_i,k_i,d_i]=[2,2,1]$ for $i=s+1$, we obtain linear codes  with parameters $[s+2, s-28, 19]$. 

\item Taking $148\leq s \leq 168$ and $[n_i,k_i,d_i]=[1,1,1]$ for $1\leq i\leq s$  and  $[n_i,k_i,d_i]=[3,2,2]$ for $i=s+1$, we obtain linear codes  with parameters $[s+3, s-28, 20]$.

\item Taking $164\leq s \leq 168$ and $[n_i,k_i,d_i]=[1,1,1]$ for $1\leq i\leq s$  and  $[n_i,k_i,d_i]=[3,2,2]$ for $i=s+1,s+2$, we obtain linear codes  with parameters $[s+6, s-26, 20]$.

\item Taking $166\leq s \leq 168$ and $[n_i,k_i,d_i]=[1,1,1]$ for $1\leq i\leq s$  and  $[n_i,k_i,d_i]=[4,2,3]$ for $i=s+1$, we obtain linear codes  with parameters $[s+4, s-28, 20]$.

\end{enumerate}

\end{example}

In example 4, we consider the function field $F=\mathbb{F}_{64}(x,y)$  defined by the equation

\[y^4+y^2+y=x^9.\] 

It can be proven that $F/\mathbb{F}_{64}$ has genus $12$, $257$ rational places and  many places of degree three.
The only pole of $x$ and $y$ will be denoted by $P_{\infty}$. 

\begin{example}
Let $A=23P_\infty$, $B=22P_\infty$, $G=A+B=45P_\infty$ and  $Z=P_{\infty}$, then $\ell(G)=34$.

\begin{enumerate}
\item Taking $228\leq s \leq 256$ and $[n_i,k_i,d_i]=[1,1,1]$ for $1\leq i\leq s$, we obtain linear codes  with parameters $[s, s-34, 24]$.

\item Taking $255\leq s \leq 256$ and $[n_i,k_i,d_i]=[1,1,1]$ for $1\leq i\leq s$  and  $[n_i,k_i,d_i]=[3,2,2]$ for $i=s+1$, we obtain linear codes  with parameters $[s+3, s-32, 24]$.

\item Taking $254\leq s \leq 256$ and $[n_i,k_i,d_i]=[1,1,1]$ for $1\leq i\leq s$  and  $[n_i,k_i,d_i]=[3,2,2]$ for $s+1 \leq i \leq s+2$, we obtain linear codes  with parameters $[s+6, s-30, 24]$.

\item Taking $253\leq s \leq 256$ and $[n_i,k_i,d_i]=[1,1,1]$ for $1\leq i\leq s$  and  $[n_i,k_i,d_i]=[3,2,2]$ for $s+1 \leq i \leq s+3$, we obtain linear codes  with parameters $[s+9, s-28, 24]$.

\item Taking $252\leq s \leq 256$ and $[n_i,k_i,d_i]=[1,1,1]$ for $1\leq i\leq s$  and  $[n_i,k_i,d_i]=[3,2,2]$ for $s+1 \leq i \leq s+4$, we obtain linear codes  with parameters $[s+12, s-26, 24]$.

\end{enumerate}

\end{example}

\bibliographystyle{plain}

\bibliography{biblio}

\end{document}